\documentclass{amsart}
\usepackage{amsmath}
\usepackage{amsthm}

\newtheorem{thm}{Theorem}[section] 
\newtheorem{pro}[thm]{Proposition}  
\newtheorem{cor}[thm]{Corollary}    

\theoremstyle{definition}           

\newtheorem{defn}[thm]{Definition}  

\numberwithin{equation}{section}
\renewcommand{\Re}{\mathop{\mathrm{Re}}\nolimits}
\newcommand{\Log}{\mathop{\mathrm{Log}}\nolimits}
\newcommand{\Arg}{\mathop{\mathrm{Arg}}\nolimits}
\renewcommand{\Im}{\mathop{\mathrm{Im}}\nolimits}
\renewcommand{\mod}{\mathop{\mathrm{mod}}\nolimits}

\newcommand{\bT}{{\mathbb T}}
\newcommand{\bD}{{\mathbb D}}
\newcommand{\bC}{{\mathbb C}}
\newcommand{\bS}{{\mathbb S}}
\newcommand{\sign}{\operatorname{sign}}
\newcommand{\bR}{{\mathbb R}}
\newcommand{\cP}{{\mathcal P}}
\newcommand{\cT}{{\mathcal T}}
\newcommand{\cG}{{\mathcal G}}
\newcommand{\cH}{{\mathcal H}}
\newcommand{\ol}{\overline}
\newcommand{\wh}{\widehat}
\renewcommand{\to}{\longrightarrow}
\renewcommand{\i}{\infty}
\renewcommand{\o}{\omega}
\renewcommand{\th}{\theta}
\renewcommand{\a}{\alpha}
\renewcommand{\b}{\beta}
\renewcommand{\l}{\lambda}
\newcommand{\var}{\varphi}

\begin{document}
\title[A Mini-Max Problem for Self-Adjoint Toeplitz Matrices]{A Mini-Max Problem for Self-Adjoint Toeplitz Matrices}

\author[Dennis Courtney and Donald Sarason]{Dennis Courtney${}^*$ and Donald Sarason}

\address{Department of Mathematics, The University of Iowa, Iowa City, IA\ \ 52242-1419, USA}

\email{dennis-courtney@uiowa.edu}

\address{Department of Mathematics, University of California, Berkeley, CA\ \ 94720-3840, USA}

\email{sarason@math.berkeley.edu}



\thanks{${}^*$Partially supported by the University of Iowa Department of Mathematics NSF VIGRE grant DMS-0602242.} 

\begin{abstract}
We study a minimum problem and associated maximum problem for finite, complex, self-adjoint Toeplitz matrices.  If $A$ is such a matrix, of size $(N+1)$-by-$(N+1)$, we identify $A$ with the operator it represents on $\cP_N$, the space of complex polynomials of degrees at most $N$, with the usual Hilbert space structure it inherits as a subspace of $L^2$ of the unit circle.  The operator $A$ is the compression to $\cP_N$ of the multiplication operator on $L^2$ induced by any function in $L^{\i}$ whose Fourier coefficients of indices between $-N$ and $N$ match the matrix entries of $A$.  Our minimum problem is to minimize the $L^{\i}$ norm of such inducers.  We show there is a unique one of minimum norm, and we describe it.  The associated maximum problem asks for the maximum of the ratio of the preceding minimum to the operator norm.  That problem remains largely open.  We present some suggestive numerical evidence.
\end{abstract}

\maketitle

\section{Introduction}
\label{sec1}

To begin we consider an $(N+1)$-by-$(N+1)$ complex Toeplitz matrix $A$ ($N$ a positive integer):
\[
A = (a_{j-k})_{j,k=0}^N.
\]
By $L^2$ we shall mean the $L^2$ space of normalized Lebesgue measure on the unit circle, $\bT$.  The subspace of $L^2$ consisting of the complex polynomials of degrees at most $N$ will be denoted by $\cP_N$.  We identify $A$ with the operator it induces on $\cP_N$, i.e., the operator on $\cP_N$ whose matrix with respect to the monomial basis is $A$.

The operator $A$ is the compression to $\cP_N$ of the multiplication operator on $L^2$ induced by any function $f$ in $L^{\i}$ of the circle whose Fourier coefficients indexed between $-N$ and $N$ match the matrix elements of $A$, in other words, for which $\wh{f}(n) = a_n$ for $n = -N,-N+1,\dots,0,\dots,N-1,N$.  The family of such functions $f$ will be denoted by $\cG_A$.  We pose the problem of minimizing $\|f\|_{\i}$ over the class $\cG_A$.  Standard reasoning shows the minimum is attained; we denote it by $c_A$.  It is evident that $\|A\| \le c_A$.

In the cases where $A$ is lower triangular or upper triangular, our problem fits into the classical interpolation problem of C.~Carath\'eodory and L.~Fej\'er, dating back 100 years \cite{CF}.  Perhaps surprisingly, the more general problem seems not to have been considered until recently \cite{BT}, \cite{NF}, \cite{S}.  Although we will not do so here, it can be treated as a special case of a general dual extremal problem studied by M. G. Krein and A. A. Nudelman in \cite{KN} --- more details later.

In this paper we focus on the case where $A$ is self-adjoint.  For this case we show that the minimizing function in $\cG_A$ is unique, and we describe it: it alternates the values $c_A$ and $-c_A$ on a family of subarcs that partition $\bT$, are even in number, and number at most $2N$.

If the matrix $A$ is self-adjoint but not diagonal, then the minimizing function is not a rational function and so multiplies every nonzero function in $\cP_N$ outside of $\cP_N$.
It follows that $c_A > \|A\|$.  Let $\cT_{N+1}$ denote the class of $(N+1)$-by-$(N+1)$ self-adjoint Toeplitz matrices.  We pose the problem of maximizing the ratio $c_A/\|A\|$ over the nonzero matrices $A$ in $\cT_{N+1}$.  As in the minimum problem, standard reasoning shows that the maximum is attained; we denote it by $c_N$.

It is proved in \cite{S} that $c_1 = \pi/2$.  The other cases remain open.  We have numerical evidence, presented in Section~\ref{sec7}, that $c_N > \pi/2$ for $N > 1$, and suggesting that $c_N$ increases with $N$.  We are able to prove that the preceding inequality holds for infinitely many $N$.

It is known, and nontrivial, that the numbers $c_N$ have a common bound.  \linebreak M.~Bakonyi and D.~Timotin obtain in \cite{BT} the bound $c_N \le 2$.  The interpolation problem they consider is slightly different from ours, but their reasoning applies in our case.  L.~N. Nikolskaya and Yu.~B. Farforovskaya obtain in \cite{NF} the bound $c_N \le 3$, and their reasoning applies even in the non-self-adjoint case.  In both papers \cite{BT} and \cite{NF}, the inequality obtained is an outcome of a matrix extension problem.

In Section~\ref{sec2} we state what we shall need concerning the Carath\'eodory--Fej\'er interpolation problem.  Section~\ref{sec3} recasts our minimum problem in a form that fits into the Carath\'eodory--Fej\'er framework.  The minimizing functions for our minimum problem are identified in Section~\ref{sec4} and studied further in Section~\ref{sec5}.  The inequality $c_1 \le \pi/2$ is reproved in Section~\ref{sec6} in a way that exploits the results from Sections~\ref{sec4} and \ref{sec5}.  The concluding Section~\ref{sec8} contains open questions and conjectures.  The maximum part of our mini-max problem remains largely mysterious, while strongly
 enticing.  New ideas are needed.

We are indebted to an anonymous referee for alerting us to the relevance to our study of references \cite{G} and \cite{KN}.

\section{Carath\'eodory--Fej\'er Interpolation Problem}
\label{sec2}

The Carath\'eodory--Fej\'er problem asks whether there is a holomorphic function in the unit disk $\bD$ having prescribed power series coefficients of orders $0,1,\dots,N$ and having a prescribed supremum norm.  Carath\'eodory and Fej\'er's solution, recast in our context and in our notation, can be stated as follows:

\medbreak
{\sc Theorem CF.} {\em 
If $A$ is lower triangular then $c_A = \|A\|$, attained by a unique function in $\cG_A$, a Blaschke product of order at most $N$ multiplied by $\|A\|$.
}

\medbreak
This is actually a mild generalization of what Carath\'eodory and Fej\'er prove, because their minimization problem is over $H^{\i}$ functions, not $L^{\i}$ functions.  But once one knows their result, one easily generalizes it to the one stated above.  Theorem~CF forms the basis of the analysis to follow.

\section{Reformulation of the Minimum Problem}
\label{sec3}

Let $A$ be a matrix in the class $\cT_{N+1}$, and let the function $f$ belong to $\cG_A$.  Then $\Re f$ is also in $\cG_A$ and $\|\Re f\|_{\i} \le \|f\|_{\i}$.  Hence $c_A$, the minimum of $\|f\|_{\i}$ over the class $\cG_A$, equals the minimum of $\|f\|_{\i}$ over the class $\cG_A^* = \{f \in \cG_A: f = \ol{f}\}$.

We let $P_+$ denote the orthogonal projection of $L^2$ onto the Hardy space $H^2 = \{h \in L^2: \wh{h}(n) = 0$ for $n < 0\}$.  As is customary, we identify the functions in $H^2$ with their holomorphic extensions in $\bD$.

Suppose $f$ is a function in the class $\cG_A^*$.  Then the function $h = P_+f - a_0/2$ is in $H^2$, and $f = h + \ol{h} = 2\Re h$.  We have
\[
h(z) = a_0/2 + a_1z + \dots + a_Nz^N + O(z^{N+1})\quad (z \to 0).
\]
Let $\cH_A$ denote the class of functions $h$ in $H^2$ that satisfy the preceding condition and for which $\Re h$ is in $L^{\i}$.  Then $2\Re h$ is in $\cG_A^*$ for such an $h$.  Accordingly, our minimum problem can be restated as the problem of minimizing $2\|\Re h\|_{\i}$ over the class $\cH_A$.

\section{The Minimum Problem}
\label{sec4}

Let $A$ be a matrix in $\cT_{N+1}$.  To eliminate a trivial case, we assume here that $A$ is not diagonal.  Suppose $h$ is a minimizing function in $\cH_A$, i.e., a function in $\cH_A$ satisfying $2\|\Re h\|_{\i} = c_A$.  To simplify the notation we temporarily write $c$ in place of $c_A$ below.

We introduce the domain
\[
\bS_c = \left\{z \in \bC: -\frac {c}{2} < \Re z < \frac {c}{2}\right\},
\]
an open vertical strip in the plane bisected by the imaginary axis, of width $c$.  The domains $\bD$ and $\bS_c$ are conformally equivalent under the map $h_c: \bD \to \bS_c$ given by
\[
h_c(z) = \frac {c}{\pi i} \Log\left( \frac {1+z}{1-z}\right);
\]
here, $\Log$ denotes the principal branch of $\log$.  The inverse map $h_c^{-1}: \bS_c \to \bD$ is given by
\[
(h_c^{-1})(z) = i \tan\left( \frac {\pi z}{2c}\right).
\]

The composite function
\[
b = h_c^{-1} \circ h = i \tan\left( \frac {\pi h}{2c}\right)
\]
is a self-map of $\bD$.  We write its power series about $0$ as
\[
b(z) = b_0 + b_1z + \dots + b_Nz^N + O(z^{N+1}).
\]

The relation between the coefficients $a_0,\dots,a_N$ of $h$ and $b_0,\dots,b_N$ of $b$ is rather complicated.  To illustrate, we exhibit below the expressions for $b_0,b_1,b_2,b_3$ in terms of $a_0,a_1,a_2,a_3$.
\[
\begin{aligned}
b_0 &= i\tan\left( \frac {\pi a_0}{4c} \right),\\
b_1 &= \frac {\pi ia_1}{2c} \sec^2 \left( \frac {\pi a_0}{4c} \right), \\
b_2 &= \frac {\pi ia_2}{2c} \sec^2 \left( \frac {\pi a_0}{4c} \right) + \frac {\pi^2ia_1}{2^2c^2} \sec^2 \left( \frac {\pi a_0}{4c} \right) \tan \left( \frac {\pi a_0}{4c} \right), \\
b_3 &= \frac {\pi ia_3}{2c} \sec^2 \left( \frac {\pi a_0}{4c} \right) + \frac {\pi^2a_1a_2}{2^2c^2} \sec^2 \left( \frac {\pi a_0}{4c} \right) \tan \left( \frac {\pi a_0}{4c} \right) \\
&\quad + \frac {\pi^3ia_1^3}{3 \cdot 2^3 \cdot c^3} \left[ 2\sec^2\left( \frac {\pi a_0}{4c} \right) \tan^2 \left( \frac {\pi a_0}{4c} \right) + \sec^4\left( \frac {\pi a_0}{4c} \right)\right].
\end{aligned}
\]
However, for our purpose these complications are not an obstacle.  From the rule for composing formal power series one knows that, for $0 \le n \le N$, the coefficient $b_n$ depends only on $c,a_0,\dots,a_n$, and vice versa.  That is all we shall need.

Let $B$ denote the $(N+1)$-by-$(N+1)$ lower triangular Toeplitz matrix with entries $b_0,b_1,\dots,b_N$ in the first column.  By Theorem~CF, there is a unique Blaschke product $\o$ of order $\le N$ such that
\[
\|B\|\o(z) = b_0 + b_1z + \dots + b_Nz^N + O(z^{N+1}).
\]
We then have
\[
((h_c^{-1}) \circ (\|B\|\o)) = \frac {a_0}{2} + a_1z + \dots + a_Nz^N + O(z^{N+1}).
\]
But in fact $\|B\| = 1$, for the inequality $\|B\| < 1$ would imply that twice the real part of the preceding composite function has $L^{\i}$-norm less than $c = c_A$, contradicting the assumed minimizing property of $h$.  Hence
\[
h = h_c \circ \o = \frac {c_A}{\pi i} \Log\left( \frac {1+\o}{1-\o}\right).
\]

\begin{defn}
\label{def4.1}
For $c > 0$ and $n$ a positive integer, an alternating step function of height $c$ and order $n$ is a function in $L^{\i}$ that assumes alternatively the values $c$ and $-c$ on $2n$ subarcs that form a partition of $\bT$.
\end{defn}

For example, the function
\[
\psi(e^{i\th}) = \begin{cases}
1, &0 < \th < \pi, \\
-1, &-\pi < \th < 0,
\end{cases}
\]
is an alternating step function of height $1$ and order $1$.

For the function $h = h_c \circ \o$ above, we have, on $\bT$ (with $\Arg$ denoting the principal value of $\arg$),
\[
\begin{aligned}
2\Re h &= \frac {2c_A}{\pi} \Arg\left( \frac {1+\o}{1-\o} \right) = \frac {2c_A}{\pi} \Arg((1+\o)(1-\ol{\o})) \\
&= \frac {2c_A}{\pi} \Arg(\o-\ol{\o}) = \frac {2c_A}{\pi} \Arg(i\Im \o) \\
&= \frac {2c_A}{\pi} \left( \frac {\pi}{2} \begin{cases}
+0 &\text{if $\sin \o > 0$,} \\
-\pi &\text{if $\sin \o < 0$,}
\end{cases} \right) \\
&= \begin{cases}
c_A &\text{if $\sin \o > 0$,} \\
-c_A &\text{if $\sin \o < 0$.}
\end{cases}
\end{aligned}
\]
Let $n$ be the order of $\o$.  The Blaschke product $\o$ wraps the unit circle $\bT$ around itself $n$ times in the strictly counterclockwise direction.  We thus see that the minimizing function $2\Re h$ is an alternating step function of height $c_A$ and order $n$.

The following theorem summarizes what has been proved.  To include the diagonal case in the statement, we define an alternating step function of order $0$ to be a constant function.  The constant is defined to be the height of the function.

\begin{thm}
\label{th4.1}
For $A$ a matrix in $\cT_{N+1}$, there is a unique function in $\cG_A$ of minimum $L^{\i}$-norm.  It is an alternating step function of height $c_A$ and order at most $N$.
\end{thm}

Theorem~\ref{th4.1} as stated excludes much of the information obtained in the reasoning leading up to it--- information we shall use later.  The theorem as stated can be obtained by fitting it into the general dual extremum problem studied by Krein and Nudelman in Chapter IX of \cite{KN}.  The general problem, which involves a Banach space and its dual, is formulated at the beginning of \S 1 of Chapter IX, after which their basic theorem on the problem, Theorem 1.1, is established.  Taking real $L^1$ of the unit circle as the Banach space in question, it follows fairly simply from this theorem that each of our minimizers has the form $L \sign \nu$, where $L$ is a positive number and $\nu$ is a real trigonometric polynomial of order at most $N$; the minimizer is thus an alternating step function.  Establishing the uniqueness of the minimizer takes more work but can be done on the basis of Krein-Nudelman's Theorem 1.3 and its corollaries.  See also Theorem 2.1 and 2.2 from Chapter IX of \cite{KN}.  Of course, the approach we use here, which  is specific to our setting, is quite distinct from the approach in \cite{KN}.

The reasoning that led to the theorem can be reversed.  Start with a nonconstant Blaschke product $\o$ of order $n \le N$.  Fix a positive number $c$, and form the function $h = \frac {c}{\pi i} \Log\left( \frac {1+\o}{1-\o}\right)$.  Then $2\Re h$ is an alternating step function of height $c$ and order $n$.  We show that $2\Re h$ is the minimizing function for the corresponding minimum problem.

In view of the analysis that produced Theorem~\ref{th4.1}, it will suffice to show that, among all $H^{\i}$ functions $g$ satisfying $\wh{g}(k) = \wh{\o}(k)$ for $k = 0,\dots,N$, the function $\o$ has the least norm.  By Theorem~CF, this amounts to showing that the compression $B$ to $\cP_N$ of the operator on $L^2$ of multiplication by $\o$ has norm $1$.

Let $z_1,\dots,z_n$ be the zeros of $\o$.  Then $\o$ has the form
\[
\o(z) = \l^2 \prod_{j=1}^n \left( \frac {z-z_j}{1-\ol{z}_jz} \right),
\]
where $\l$ is a constant of modulus~$1$.  Defining the polynomials $r$ and $s$ by
\[
r(z) = \l \prod_{j=1}^n (z-z_j),\quad s(z) = \ol{\l} \prod_{j=1}^n (1-\ol{z}_jz),
\]
we have $\o = r/s$.  The polynomials $r$ and $s$ belong to $\cP_N$, and we have $\|r\|_2 = \|s\|_2$, plus $Bs = r$.  The desired equality, $\|B\| = 1$, follows.  (The preceding argument is of course well known.)

The question whether every nonconstant alternating step function is obtained by the preceding procedure now arises.  That question is addressed in the next section.

\section{Inner Functions and Alternating Step Functions}
\label{sec5}

\begin{thm}
\label{th5.1} 
Let $\psi$ be an alternating step function of height $c$ and order $N > 0$.  Then there is a Blaschke product $\o$ of order $N$ such that
\[
\psi = \frac {2c}{\pi} \Arg\left( \frac {1+\o}{1-\o}\right)
\]
on $\bT$.
\end{thm}

\begin{proof}
We assume (without loss of generality) that $c = 1$.
\renewcommand{\qedsymbol}{}\end{proof}

\medbreak
{\sc Step 1.} Let the real numbers $\a_1,\dots,\a_N,\b_1,\dots,\b_N$ satisfy
\[
0 \le \a_1 < \b_1 < \a_2 < \b_2 < \dots < \a_N < \b_N < \a_1 + 2\pi := \a_{N+1}.
\]
Let the arcs $\chi_n^-$ and $\chi_n^+$ $(n = 1,\dots,N)$ be defined by
\[
\chi_n^- = \{e^{i\th}: \a_n < \th < \b_n\},\quad \chi_n^+ = \{e^{i\th}: \b_n < \th < \a_{n+1}\}.
\]
The function $\psi$ that takes the value $1$ on each arc $\chi_n^-$ and the value $-1$ on each arc $\chi_n^+$ is, to within a sign, the general alternating step function of height~$1$ and order~$N$.  We let
\[
\a = \a_1 + \dots + \a_N,\quad \b = \b_1 + \dots + \b_N.
\]

\medbreak
{\sc Step 2.} We introduce the polynomials
\[
p(z) = e^{-i\a/2} \prod_{n=1}^N (z-e^{i\a_n}),\quad q(z) = e^{-i\b/2} \prod_{n=1}^N (z-e^{i\b_n}).
\]
It is asserted that the polynomial $q - ip$ has degree $N$.  In fact, the leading coefficient of $q - ip$ is $e^{-i\b/2} - ie^{-i\a/2}$.  For this to vanish we must have $e^{i(\b-\a)/2} = -i = e^{3\pi i/2}$.  However $0 < \b - \a < 2\pi$, so $0 < \frac {\b-\a}{2} < \pi$, and the assertion follows.

\medbreak
{\sc Step 3.} We show that the rational function $q/p$ is real valued on $\bT$ (except at its poles $e^{i\a_n}$, $n = 1,\dots,N$), and that on each arc $\chi_n := \chi_n^- \cup \chi_n^+ \cup \{e^{i\b_n}\}$ it is an increasing function with range $\bR$.

We have
\[
\begin{aligned}
q(e^{i\th})/p(e^{i\th}) &= e^{i(\b-\a)/2} \prod_{n=1}^N (e^{i\th}-e^{i\b_n})\left/ \prod_{n=1}^N (e^{i\th}-e^{i\a_n}) \right. \\
&= e^{i(\b-\a)/2} \prod_{n=1}^N e^{i\left( \frac {\a_n-\b_n}{2}\right)} \left[ \frac {e^{i\left( \frac {\th-\b_n}{2}\right)} - e^{-i\left( \frac {\th-\b_n}{2} \right)}}{e^{i\left( \frac {\th-\a_n}{2} \right)} - e^{-i\left( \frac {\th-\a_n}{2}\right)}} \right] \\
&= \prod_{n=1}^N \frac {\sin\left( \frac {\th-\b_n}{2}\right)}{\sin \left( \frac {\th-\a_n}{2}\right)}\,.
\end{aligned}
\]
Thus $q/p$ is real valued on $\bT$.

For $\a_n < \th < \a_{n+1}$, the function $\sin\left( \frac {\th-\a_n}{2}\right)$ is positive, while the function $\sin\left( \frac {\th-\b_n}{2}\right)$ is negative for $\th < \b_n$ and positive for $\th > \b_n$.  For $j \ne n$ and $\a_n < \th < \a_{n+1}$, the functions $\sin\left( \frac {\th-\a_j}{2}\right)$ and $\sin\left( \frac {\th-\b_j}{2}\right)$ have the same sign (negative for $j < n$, positive for $j > n$).  Together with the expression above for $q/p$, this tells us that $q/p$ is negative on $\chi_n^-$ and positive on $\chi_n^+$.  As $e^{i\a_n}$ and $e^{i\a_{n+1}}$ are poles of $q/p$, we thus must have $q(e^{i\th})/p(e^{i\th}) \to -\i$ as $\th \searrow \a_n$ and $q(e^{i\th})/p(e^{i\th}) \to +\i$ as $\th \nearrow \a_{n+1}$.  Hence $q/p$ maps the arc $\chi_n$ onto $\bR$.  It must be monotone on each such arc because, being a rational function of degree $N$, it assumes no value with multiplicity greater than $N$.

\medbreak
{\sc Step 4.} We show that the function $\o = (q-ip)/(q+ip)$ is a Blaschke product of order $N$ satisfying
\[
\psi = \frac {2}{\pi} \Arg \left( \frac {1+\o}{1-\o}\right). \leqno(1)
\]

In fact, since $\o$ can be rewritten as $\o = \left( \frac {q}{p} - i \right)\left/\left( \frac {q}{p} + i\right)\right.$, it follows by Step~3 that $\o$ is unimodular on $\bT$.  From the definition of $\o$ we see immediately that $\o(e^{i\a_n}) = 1$ for each $n$ and $\o(e^{i\b_n}) = -1$ for each $n$.

We can re-express $\o$ as $\var \circ \left( \frac {q}{p}\right)$, where $\var(z) = (z-i)/(z+i)$, a linear-fractional transformation that maps the upper half-plane to $\bD$ and $\bR$ to $\bT\setminus \{1\}$, with $\var(z)$ moving counterclockwise on $\bT$ as $z$ moves in the positive direction on $\bR$.  We know from Step~3 that on each arc $\chi_n$ the function $q/p$ is increasing with range $\bR$.  We can conclude that $\arg \o$ undergoes an increment of $2\pi$ on each arc $\chi_n$, hence an increment of $2\pi N$ on $\bT$.  By the argument principle, $\o$ has $N$ zeros in $\bD$, thus is a Blaschke product of order $N$.

Because $\frac {1+\o}{1-\o} = \frac {q}{ip}$, if $q(e^{i\th})/p(e^{i\th}) > 0$, i.e., if (according to Step~3) $e^{i\th}$ is in one of the arcs $\chi_n^+$, we have $\Arg\left( \frac {1+\o(e^{i\th})}{1-\o(e^{i\th})}\right) = -\frac {\pi}{2} + 0 = -\frac {\pi}{2}$.  Similarly, if $e^{i\th}$ is in one of the arcs $\chi_n^-$ then $\Arg\left( \frac {1+\o(e^{i\th})}{1-\o(e^{i\th})}\right) = -\frac {\pi}{2} + \pi = \frac {\pi}{2}$.  The desired equality $(1)$ follows.  This concludes the proof of Theorem~\ref{th5.1}.  \qed

\medbreak

Theorem~\ref{th5.1} is a particular case of Theorem 9 in the paper \cite{G} of P. Gorkin and R.C. Rhoades.  Both \cite[Theorem 9]{G} and our Theorem~\ref{th5.1} assert the extistence of certain finite Blaschke products.  The proofs of the two theorems are essentially the same, the chief difference being that Gorkin-Rhoades work mainly on the real line while we stick to the unit circle.  A conformal isomorphism between the unit disk and the upper half-plane connects the two approaches.

We have chosen to include the proof above of Theorem~\ref{th5.1}, rather than refer the reader to \cite{G} for a proof, for the sake of completeness, and because the aims of \cite{G} are rather different from ours.  (Since Theorem~\ref{th5.1} is less general than \cite[Theorem 9]{G}, its proof is also a bit simpler.)

To conclude this section we discuss the extent to which the Blaschke product $\o$ produced in the preceding proof fails to be unique.  Note that $\left| \frac {1+\o(0)}{1-\o(0)}\right| = 1$.  In fact
\[
\frac {1+\o(0)}{1-\o(0)} = \frac {q(0)}{ip(0)} = -ie^{i(\b-\a)/2}.
\]

Suppose the Blaschke product $\var$ of order $n$ satisfies $(1)$ when substituted for $\o$.  Then the argument of $\frac {1+\var}{1-\var}$ must undergo a jump of $-\pi$ across each point $e^{i\a_n}$ and a jump of $\pi$ across each point $e^{i\b_n}$.  This implies, by a standard argument, that each point $e^{i\a_n}$ is a simple pole of $\frac {1+\var}{1-\var}$ and each point $e^{i\b_n}$ is a simple zero of $\frac {1+\var}{1-\var}$.  The functions $\frac {1+\var}{1-\var}$ and $\frac {1+\o}{1-\o}$ are thus rational functions with the same zeros and poles, so they are constant multiples of each other.  It follows that one can obtain $\var$ by composing $\o$ from the left with a conformal automorphism of $\bD$.  That conformal automorphism must fix the points $1$ and $-1$, so it has the form $z \mapsto \frac {z+r}{1+rz}$, where $-1 < r < 1$.  We get $\var = \frac {\o+r}{1+r\o}$.

A calculation produces the equality
\[
\frac {1+\var}{1-\var} = \left( \frac {1+r}{1-r} \right) \left( \frac {1+\o}{1-\o}\right).
\]
Thus $\left| \frac {1+\var(0)}{1-\var(0)}\right| = 1$ only if $r = 0$, i.e., only if $\var = \o$.  We see that $\o$ becomes unique if one imposes the supplementary condition $\left| \frac {1+\o(0)}{1-\o(0)} \right| = 1$.

\begin{thm} 
\label{th5.2}
If $\psi$ is an alternating step function of order $N$, then $\psi$ is uniquely determined by the Fourier coefficients $\wh{\psi}(0),\wh{\psi}(1),\dots,\wh{\psi}(N)$, in the sense that if those coefficients are known then the coefficients $\wh{\psi}(N+1),\wh{\psi}(N+2),\dots$ can in principle be found.
\end{thm}

\begin{proof}
Let $\psi$ be such a function, and let $c$ be its height.  We may assume $N > 0$ (the case $N = 0$ being trivial).  By Theorem~\ref{th5.1}, there is a Blaschke product $\o$ of order $N$ such that $\psi = h + \ol{h}$, where 
\[
h = \frac {c}{\pi i} \Log\left( \frac {1+\o}{1-\o} \right).
\]
From this we see that, for any $n$, the coefficient $\wh{\psi}(n)$ can be expressed in terms of $c$ and the coefficients $\wh{\o}(0),\dots,\wh{\o}(n)$.  As $\o$ is a rational function of degree $N$, the coefficients $\wh{\o}(N+1),\wh{\o}(N+2),\dots$ can be expressed in terms of $\wh{\o}(0),\dots,\wh{\o}(N)$.  The desired conclusion follows.
\end{proof}

\begin{cor}
\label{cor5.3}
Let $E$ be the union of $N$ closed disjoint subarcs of $\bT$.  Then $E$ is uniquely determined by the coefficients $\wh{\chi}_E(0),\dots,\wh{\chi}_E(N)$.
\end{cor}

\begin{proof}
It suffices to apply Theorem~\ref{th5.2} to the alternating step function $2\chi_E-1$.
\end{proof}

The phenomenon described in the corollary is treated at length in \cite{C}.

\section{Two Dimensions}
\label{sec6}

We indicate how the equality $c_1 = \frac {\pi}{2}$ can be deduced from our characterization of minimizing functions.  The inequality states that $c_A/\|A\| \le \frac {\pi}{2}$ for all nonzero matrices $A$ in $\cT_2$.  For such a matrix $A$ we know from Theorem~\ref{th4.1} that the minimizing function in $\cG_A$ is an alternating step function of order $1$ and height $c_A$.  We can thus try to establish the inequality by starting with an alternating step function $\psi$ of order $1$ and a given height, computing the corresponding matrix $A$, finding $\|A\|$, and checking the ratio of interest.  And it will clearly suffice to treat the case where $\psi$ has height $1$ and has a discontinuity at the point $1$ with a jump there of $2$.

Let $\a$ be a point in the interval $(0,2\pi)$, and define $\psi$ by
\[
\psi(e^{i\th}) = \begin{cases}
1, &0 < \th < \a, \\
-1, &\a < \th < 2\pi.
\end{cases}
\]
Then $\psi$ is the general alternating step function of interest.  The corresponding matrix $A$ is given by
\[
A = \begin{pmatrix}
\wh{\psi}(0) & \wh{\psi}(-1) \\
\wh{\psi}(1) & \wh{\psi}(0)
\end{pmatrix}\,.
\]
Calculations give
\[
\wh{\psi}(0) = \frac {\a-\pi}{\pi},\quad \wh{\psi}(1) = \frac {2e^{-i\a/2}\sin \frac {\a}{2}}{\pi}\,.
\]
Thus
\[
A = \frac {1}{\pi} \begin{pmatrix}
\a-\pi & 2e^{i\a/2}\sin \frac {\a}{2} \\
2e^{-i\a/2}\sin \frac {\a}{2} & \a - \pi
\end{pmatrix}\,.
\]
The eigenvalues of $A$ can be found by the standard procedure.  The calculations yield
\[
\|A\| = \frac {|\a-\pi| + 2\left| \sin \frac {\a}{2}\right|}{\pi}\,.
\]
As $c_A = 1$, the desired conclusion is $\frac {1}{\|A\|} \le \frac {\pi}{2}$.  We have
\[
\frac {1}{\|A\|} = \frac {\pi}{|\a-\pi| + 2\left| \sin \frac {\a}{2}\right|}\,,
\]
so we need only verify that
\[
\pi - \a + 2\sin \frac {\a}{2} \ge 2
\]
for $0 \le \a \le \pi$.  Differentiation with respect to $\a$ shows that the function on the left side is decreasing on $[0,\pi]$.  At $\a = \pi$ it has the value $2$, so the desired inequality holds, and it reduces to an equality only for $\a = \pi$.

\section{The Numbers $c_N$}
\label{sec7}

We present here what we have learned about the maxima $c_N$, including numerical results for small $N > 1$.  For $f$ a function in $L^{\i}$ we let $A_{f,N}$ denote the compression to $\cP_N$ of the multiplication operator on $L^2$ induced by $f$.

The inequality $c_N \ge \frac {\pi}{2}$ is easily seen at this point.  For example, let $\psi$ be the alternating step function of height $1$ and order $N$ associated with the inner function $\o(z) = z^N$ as in the proof of Theorem~\ref{th4.1}:
\[
\psi(z) = 2\Re\left( \frac {1}{\pi i} \Log \left( \frac {1+z^N}{1-z^N}\right)\right).
\]
Referring to that proof, one sees that $\psi$ alternates the values $1$ and $-1$ on $2N$ consecutive subarcs of $\bT$, each of length $\frac {\pi}{N}$.  For $z \to 0$ we have
\[
\frac {1+z^N}{1-z^N} = 1 + 2z^N + O(|z|^{2N}),
\]
so
\[
\begin{aligned}
\psi(z) &= \frac {2}{\pi} \Im\left( \Log\left( \frac {1+z^N}{1-z^N}\right)\right) \\
&= \frac {2}{\pi} \Im (\Log(1+2z^N)) + O(|z|^{2N}) \\
&= \frac {2}{\pi} \Im (2z^N) + O(|z|^{2N}) \\
&= \frac {2}{\pi} \left( \frac {2z^N-2\ol{z}^N}{2i}\right) + O(|z|^{2N}) \\
&= \frac {2}{\pi i}(z^N - \ol{z}^N) + O(|z|^{2N}).
\end{aligned}
\]
We see that the matrix $A_{\psi,N}$ has only two nonzero entries, $\frac {2}{\pi i}$ in the lower left corner and $-\frac {2}{\pi i}$ in the upper right corner.  In particular, $\|A_{\psi}\| = \frac {2}{\pi}$, and the ratio $\frac {\|\psi\|_{\i}}{\|A_{\psi,N}\|} = \frac {1}{\|A_{\psi,N}\|}$ equals $\frac {\pi}{2}$, telling us that $c_N \ge \frac {\pi}{2}$.

The following observation produces a generalization of the last inequality.

\begin{pro}
\label{prop7.1}
Let the function $f$ be in $L^{\i}$, let $k$ be a positive integer, and let the function $g$ on $\bT$ be defined by $g(e^{i\th}) = f(e^{ik\th})$.  Then $\|A_{f,N}\| = \|A_{g,kN}\|$.
\end{pro}

\begin{proof}
To simplify the notation, let $T = A_{g,kN}$.  For $p$ in $\cP_{kN}$, the $m$-th Fourier coefficient of $Tp$ is given by
\[
(Tp)^{\wedge}(m) = \sum_{\substack{0 \le j \le kN \\ -kN \le m-j \le kN}} \wh{g}(m-j)\wh{p}(j).
\]
If this coefficient is nonzero, there must be a $j$ such that $k$ divides $m-j$, in other words, such that $j \equiv m\ (\mod k)$.

For $l = 0,\dots,k-1$, let $S_l$ be the subspace of $\cP_{kN}$ consisting of the polynomials $p$ such that $\wh{p}(j) = 0$ if $j \not\equiv l\ (\mod k)$.  These subspaces are mutually orthogonal, they span $\cP_{kN}$, and, by the observation above, they are $T$-invariant.  Hence $T$ is the direct sum of its restrictions to the subspaces $S_l$.

If $p$ is in $S_0$ then there is a $q$ in $\cP_N$ such that $p(z) = q(z^k)$.  For the $km$-th coefficient of $Tp$ we have
\[
\begin{aligned}
(Tp)^{\wedge}(km) &= \sum_{j=0}^N \wh{g}(km-kj)\wh{p}(kj) \\
&= \sum_{j=0}^N \wh{f}(m-j)\wh{q}(j),
\end{aligned}
\]
showing that $T\mid S_0$ is unitarily equivalent to $A_{f,N}$.

For $p$ in $S_l$ with $l \ne 0$, the $(l+km)$-th coefficient of $Tp$ $(m = 0,\dots,N-1)$ is given by
\[
(Tp)^{\wedge}(l+km) = \sum_{\substack{j \equiv l\ (\mod k) \\ l \le j \le l+k(N-1)}} \wh{g}(l+km-j)\wh{p}(j).
\]
In the summation on the right, $j-l$ runs through the numbers $0,k,\dots,(N-1)k$, and correspondingly, $j$ runs through the numbers $l,l+k,\dots,l+k(N-1)$.  The equality can be rewritten as
\[
(Tp)^{\wedge}(l+km) = \sum_{j=0}^{N-1} \wh{g}((m-j)k)\wh{p}(l+jk).
\]
A polynomial $p$ in $S_l$ $(l \ne 0)$ can be written as $p(z) = z^lq(z^k)$, where $q$ is in $\cP_N$ and $\wh{q}(N) = 0$.  This gives us
\[
(Tp)^{\wedge}(l+km) = \sum_{j=0}^{N-1} \wh{f}(m-j)\wh{q}(j),
\]
showing that the restriction of $T$ to $S_l$ is unitarily equivalent to the $N$-by-$N$ principal minor of $A_{f,N}$, hence is of norm at most that of $A_{f,N}$.  The equality $\|A_{f,N}\| = \|A_{g,kN}\|$ now follows.
\end{proof}

The following corollary generalizes the inequality $c_N \ge \frac {\pi}{2}$.

\begin{cor}
\label{cor7.2}
For $k$ a positive integer, $c_{kN} \ge c_N$.
\end{cor}

\begin{proof}
It suffices to apply Proposition~\ref{prop7.1} with $f$ equal to an alternating step function $\psi$ of order at most $N$ for which $\frac {\|\psi\|_{\i}}{\|A_{\psi,N}\|} = c_N$.
\end{proof}

One way to obtain better lower bounds than $c_N \ge \frac {\pi}{2}$ for $N > 1$ is to find alternating step functions $\psi$ of height $1$ and order $N$ for which $\|A_{\psi, N}\|$ can be shown to be smaller than $\frac {2}{\pi}$.  We do this for a number of values of $N$.  Our method is to consider a family of step functions involving a number of parameters, and to observe that when the parameters are judiciously chosen, the Fourier coefficients of $\psi$ are related to one another in a way that permits calculation of $\|A_{\psi,N}\|$.  A computer was used to assist in the search for relations between Fourier coefficients that would yield nontrivial lower bounds on $c_N$.

\medbreak
{\sc The case $N=2$.}  For any positive $L$ and $M$ satisfying $0 < L + M < \pi$, there is a unique alternating step function $\psi$ of order $2$ taking the value $1$ at $1$ and having discontinuities at $e^{\pm iL}$ and $e^{i(\pi \pm M)}$.  Calculations show that
\[
\wh{\psi}(0) = \frac{2L+2M-\pi}{\pi}, \quad \wh{\psi}(1) = \frac{2 (\sin L - \sin M)}{\pi}, \quad \wh{\psi}(2) = \frac{\sin(2M)+\sin(2L)}{\pi}.
\]
The values of $L$ and $M$ can additionally be chosen so that the relations $\wh{\psi}(1) = 2 \wh{\psi}(0)$ and $\wh{\psi}(2) = -\wh{\psi}(1)$ are also satisfied, and these choices are unique subject to the conditions $L \in [0,1]$ and $M \in [0,2]$.  (To abbreviate the argument: one can show that subject to the inequalities on $L$ and $M$, each equation of Fourier coefficients defines $M$ implicitly as a monotone function of $L$ on $[0,1]$, and use the intermediate value theorem to see that the functions intersect exactly once.)  These relations allow one to explicitly compute $\|A_{\psi, 2}\| = 3 |\wh{\psi}(0)|$.

\medbreak
{\sc Conclusion.} If $L$ and $M$ are the unique real numbers satisfying the relations $0 < L + M < \pi$, $L \in [0,1]$, $M \in [0,2]$, $2 L + 2M - \pi = \sin L - \sin M$, and $-2 (\sin L - \sin M) = \sin(2 M) + \sin(2L)$, then
\begin{equation}\label{c2est}
c_2 \ge \frac{\pi}{3 |2L+2M-\pi|}.
\end{equation}

\medbreak
One must find $L \approx .2138$ and $M \approx 1.0263$ to several digits to deduce that this estimate improves on the known $c_2 \ge \frac{\pi}{2}$, but it can be shown that the right hand side of \eqref{c2est} is between $1.6185$ and $1.6186$.

\medbreak
{\sc The case $N=3$.}  For any $L$ in $(0, \frac{\pi}{2})$ there is a unique alternating step function $\psi(e^{i\theta})$ of order $3$ and height $1$ taking the value $1$ at $1$ and having discontinuities at $e^{\pm iL}$, $e^{i(\pi \pm L)}$, and $e^{\pm i\pi/2}$.   Calculations show that $\wh{\psi}(0)=\wh{\psi}(2)=0$ and
\[
\wh{\psi}(1) = \frac{2}{\pi}(2 \sin L - 1), \quad \wh{\psi}(3) = \frac{2}{3\pi}(2 \sin(3L) + 1).
\]
There is a unique $L$ in $(0, \frac{\pi}{2})$ satisfying $\sin(3L)+3\sin L -1=0$.  For this value of $L$ we have $\wh{\psi}(1)=-\wh{\psi}(3)$, and therefore $\|A_{\psi,3}\| = \sqrt{2} |\wh{\psi}(1)|$.  We deduce
\[
c_3 \ge \frac{\pi}{2\sqrt{2}(1 - 2 \sin L)}.
\]
Trigonometric identities imply that $\sin L$ is the least positive root of the polynomial $1 - 6x+4x^3$.  Using this polynomial one can find a similar characterization of $\frac{1}{1-2\sin L}$.

\medbreak
{\sc Conclusion.} If $k$ is the largest root of the polynomial $p(x) = 1-3x-3x^2+3x^3$, then
\begin{equation}\label{c3est}
c_3 \ge \frac{\pi}{2 \sqrt{2}} k.
\end{equation}

\medbreak
It can be checked by hand that $p(\frac{3}{2}) < 0$, and thus $k \ge \frac{3}{2}$, so $c_3 \ge \frac{3\pi}{4 \sqrt{2}} > \frac{\pi}{2}$ can be shown by hand.  Working to greater precision one can show that the right hand side of \eqref{c3est} is between $1.6825$ and $1.6826$.

Our higher dimensional examples are similar to those just given.  The difference is that, with the exception of the case $N=5$, the number of parameters involved makes exact formulas difficult to prove.  We can nevertheless give numerical values of parameters that produce a $\psi$ with the property that $\|A_{\psi,N}\|$ is small enough to give new information about $c_N$.

\medbreak
{\sc The case $N=4$.} For any positive $L,M,N,O$ satisfying $L+M+N+O<\pi$ there is an alterating step function of height $1$ and order $4$ taking the value $1$ at $1$ and having discontinuities at $e^{iL}$, $e^{i(L+M)}$, $e^{i(L+M+N)}$, $e^{i(L+M+N+O)}$ and their complex conjugates.  Taking $L \approx .1396$, $M \approx 1.1143$, $N \approx 1.096$, and $O \approx .2724$ one obtains a function $\psi$ for which $\|A_{\psi,4}\|$ can be shown to be between $1.7065$ and $1.7066$.  These values of $L,M,N,O$ appear to approximate a tuple with the property that the corresponding list of Fourier coefficients $(\wh{\psi}(j))_{j=0}^4$ is a scalar multiple of $(1, -\sqrt{5}-1, 1-\sqrt{5},\sqrt{5}-1,1+\sqrt{5})$.  If $L,M,N,O$ can be chosen so that this holds, we would have $c_4 \ge (5 |\psi(0)|)^{-1} = \frac{\pi}{5(\pi-2M-2O)}$.

\medbreak
{\sc The case $N=5$.} For any positive $L$ and $M$ satisfying $L+M < \frac{\pi}{2}$ there is a unique alternating step function $\psi$ of height $1$ and order $5$ taking the value $1$ at $1$ and having discontinuities at $e^{iL}$, $e^{i(L+M)}$, $e^{i \pi/2}$, $e^{i(\pi - L)}$, $e^{i(\pi - L-M)}$, and their complex conjugates.  We have $\wh{\psi}(0) = \wh{\psi}(2)  = \wh{\psi}(4) = 0$ for all such $L$ and $M$, and it can be shown that there are unique $L$ and $M$ in the interval $(0,1)$ for which
\begin{equation}\label{fiverel}
\wh{\psi}(1) = -2 \wh{\psi}(3) = \wh{\psi}(5).
\end{equation}
For these values of $L$ and $M$ we deduce that $\|A_{\psi, 5}\| = \frac{3}{2} |\wh{\psi}(1)|$.

\medbreak
{\sc Conclusion.} If $L$ and $M$ are the unique numbers in $(0,1)$ for which the equations \eqref{fiverel} hold, then
\begin{equation}\label{c5est}
c_5 \ge \frac{\pi}{3 (1 +2 \sin L - 2 \sin(L+M))}.
\end{equation}

\medbreak
Numerical calculations can show that $L \approx .4304$, $M \approx .2326$, and that the right hand side of \eqref{c5est} is between $1.7353$ and $1.7354$.

\medbreak
{\sc The case $N=6$.} For any positive $L,M,N,O,P,Q$ whose sum is less than $\pi$, there is a unique alternating step function $\psi$ of height $1$ and order $6$ taking the value $1$ at $1$ and having discontinuities at the points $e^{iL}$, $e^{i(L+M)}$, $e^{i(L+M+N)}$, $e^{i(L+M+N+O)}$, $e^{i(L+M+N+O+P)}$, $e^{i(L+M+N+O+P+Q)}$ and their complex conjugates.  Taking $(L,M,N,O,P,Q) \approx (.0989,.7269,.2002,.7702,.7755,.2109)$ one obtains a function $\psi$ for which $\|A_{\psi,6}\|$ can be shown to be between $1.7504$ and $1.7505$.  These parameter values appear to approximate a tuple with the property that the list of Fourier coefficients of the corresponding function $\psi$ is proportional to $(1,a,-b,c,-c,b,-a,1)$, where $a<b<c$ are the roots of $8+4x-4x^2-x^3$.

\medbreak
{\sc The case $N=7$.} For any positive $L,M,N$ satisfying $L+M+N < \frac{\pi}{2}$, there is a unique alternating step function $\psi$ of height $1$ and order $7$ taking the value $1$ at $1$ and having discontinuities at $e^{iL}$, $e^{i(L+M)}$, $e^{i(\pi/2-N)}$, $e^{i\pi/2}$, $e^{i(\pi/2+N)}$, $e^{i(\pi-L-M)}$, $e^{i(\pi-L)}$, and their complex conjugates.   Taking $L \approx .0877$, $M \approx .6343$, and $N \approx .6713$ one obtains a function $\psi$ for which it can be shown that $\|A_{\psi,7}\|$ is between $1.7677$ and $1.7678$.  These values of $L,M,N$ appear to approximate a tuple with the property that the corresponding list of Fourier coefficients $(\wh{\psi}(j))_{j=0}^7$ is a scalar multiple of $(0,-1-\sqrt{2},0,1,0,-1,0,1+\sqrt{2})$.  If $L,M,N$ can be chosen so that this holds exactly, we can conclude that $c_7 \ge \frac{(\sqrt{2}+1)\pi}{4 \sqrt{2+\sqrt{2}} (1-2\cos N-2\sin L+2\sin(L+M))}$ for these values of $L,M,N$.

\medbreak
In summary, along with the exact value $c_1 = \frac {\pi}{2} \approx 1.5707963$, we have the following lower bounds for $c_2 - c_7$, found numerically:
\[
\begin{aligned}
c_2 &\ge 1.6185 \\
c_3 &\ge 1.6825 \\
c_4 &\ge 1.7065 \\
c_5 &\ge 1.7354 \\
c_6 &\ge 1.7505 \\
c_7 &\ge 1.7677.
\end{aligned}
\]

\section{Questions and Conjectures}
\label{sec8}

We present some problems for future research suggested by the numerical results of Section~\ref{sec7}, which tell us, in particular, that $c_N > \frac{\pi}{2}$ for $N=2,3,4,5,6,7$.  (As noted in Section~\ref{sec7}, for $N = 3$ the inequality can be verified by hand.)  By the corollary to Proposition~\ref{prop7.1}, we can conclude that $c_N > \frac{\pi}{2}$ for infinitely many $N$. We conjecture that the following question has an affirmative answer.

\medbreak
{\sc Question 1.} {\em 
Is $c_N > \frac{\pi}{2}$ for all $N > 1$?
}

\medbreak
Our estimates for $c_1$--$c_7$ suggest the following question.

\medbreak
{\sc Question 2.} {\em 
Does $c_N$ increase strictly as $N$ increases?
}

\medbreak
Again, we conjecture the answer is affirmative.

As noted in Section~\ref{sec1}, the inequality $c_N \le 2$ holds for all $N$, suggesting the next question.

\medbreak
{\sc Question 3.} {\em 
What is $\sup \{c_N: N \ge 1\}$?
}

\medbreak
We have no conjecture to offer on this one.

We have seen that the maximum ratio $c_N$ is attained by an alternating step function of order at most $N$, i.e., there is such a function $\psi$ for which $\|\psi\|_{\i}/\|A_{\psi,N}\| = c_N$.  The function $\psi$ here is not unique, because one gets the same ratio if one replaces $\psi$ by a nonzero real scalar multiple of itself, or by a composite (from the right) with a rotation of the circle, or by a composite with complex conjugation on the circle.  Two questions are suggested, which we combine into one.

\medbreak
{\sc Question 4.} {\em 
Are the alternating step functions of order at most $N$ that attain the maximum ratio $c_N$ actually of order $N$, and are any two of them related via the transformations listed above?
}

\medbreak
Otherwise put, the question asks whether an alternating step function attaining the ratio $c_N$ has $2N$ discontinuities, and whether the set of discontinuities is unique to within rotations and reflections on $\bT$.  We conjecture the answer is affirmative.

The questions above are of course part of the underlying problem of ``finding,'' in some reasonable sense, the maxima $c_N$ and the corresponding maximizers.  There must be an interesting underlying structure, one would think, which so far escapes us.  Answers to one or more of Questions~1--4 should provide clues to that structure.

\end{document}